\documentclass[11pt]{amsart}

\usepackage{amsmath,amsfonts,amscd,amssymb}
\theoremstyle{plain}
\newcounter{thmcount}
\newtheorem{theorem}[thmcount]{Theorem}
\newtheorem{conjecture}[thmcount]{Conjecture}

\newtheorem{lemma}[thmcount]{Lemma}
\newtheorem{corollary}[thmcount]{Corollary}
\theoremstyle{definition}

\newtheorem*{notation}{Notation}

\def\O{{\mathcal O}}
\def\cF{{\mathcal F}}
\def\F{{\mathbb F}}

\def\Q{{\mathbb Q}}

\def\R{{\mathbb R}}
\def\C{{\mathbb C}}
\def\Z{{\mathbb Z}}

\def\m{{\mathfrak m_F}}
\def\mfnr{{\mathfrak m_{F^{nr}}}}

\def\smallmatrix#1#2#3#4{{\scriptstyle\Bigl(\begin{matrix}%
  \scriptstyle #1&\scriptstyle #2\cr\scriptstyle #3&\scriptstyle #4\cr
  \end{matrix}\Bigr)}}
\def\smallmatrix#1#2#3#4{
  \genfrac{(}{.}{0pt}{1}{#1}{#3}
  \genfrac{.}{)}{0pt}{1}{#2}{#4}
}

\def\f#1{{#1}_\phi}
\def\fKv{K_{v,\phi}}

\input cyracc.def       
\font\tencyr=wncyr10
\def\sha{\text{\tencyr\cyracc{Sh}}}

\def\newmathop#1{\expandafter\gdef\csname #1\endcsname{\mathop{\rm #1}\nolimits}}
\newmathop{tr}
\newmathop{corank}
\newmathop{rk}
\newmathop{ord}
\newmathop{Gal}
\newmathop{coker}
\newmathop{Sel}
\newmathop{Ind}
\newmathop{Frob}
\newmathop{GL}
\newmathop{SL}
\newmathop{Aut}

\def\compose{{_{^{\text{o}}}}}

\let\lar\longrightarrow
\let\iso\cong
\let\tensor\otimes

\def\notdiv{\hbox{$\not|\,$}}
\def\1{{\mathbf 1}}

\def\beq#1{\begin{equation}\label{#1}}
\def\eeq{\end{equation}}
\def\beql#1{\begin{equation}\label{#1}\begin{array}{llllll}}
\def\eeql{\end{array}\end{equation}}

\catcode`\@=11
\newwrite\@corrfile
\openout\@corrfile=correct.inc
\def\corr#1#2{\marginpar{#1}%
   \write\@corrfile{\upsilon{#1}{#2}}%
}
\catcode`\@=12

\def\corr#1#2{}

\begin{document}

\title{Parity of ranks for elliptic curves with a cyclic isogeny}
\author{Tim and Vladimir Dokchitser}
\address{Robinson College, Cambridge CB3 9AN, United Kingdom}
\email{t.dokchitser@dpmms.cam.ac.uk}
\address{Max Planck Institut f\"ur Mathematik, D-53111 Bonn, Germany}
\email{v.dokchitser@dpmms.cam.ac.uk}

\begin{abstract}
Let $E$ be an elliptic curve over a number field $K$ which admits a cyclic
$p$-isogeny with $p\ge 3$ and semistable
at primes above $p$. We determine the root number and the
parity of the $p$-Selmer rank for $E/K$, in particular confirming the parity
conjecture for such curves. We prove the analogous results for $p=2$
under the additional assumption that $E$ is not supersingular at primes
above 2.
\end{abstract}


\maketitle

\section{Introduction}

If $E$ is an elliptic curve over a number field $K$, the number of copies
of $\Z$ in the group of rational points $E(K)$ is called the Mordell-Weil rank
of $E/K$. If the Tate-Shafarevich group $\sha(E/K)$ is finite
(conjecturally, this is always the case), then for every prime $p$
it is the same as the $p$-Selmer rank of $E/K$, defined
as the Mordell-Weil rank
plus the number of copies of $\Q_p/\Z_p$ in $\sha(E/K)$.
We will be concerned with the parity of the $p$-Selmer rank,
and we will write $\sigma(E/K,p)$ for $(-1)^{\text{$p$-Selmer rank of $E/K$}}$.

Tate's generalization of the Birch and Swinnerton-Dyer conjecture for
elliptic curves over number fields predicts that the Mordell-Weil rank
is the same as the analytic rank, the order of vanishing of the
$L$-function $L(E/K,s)$ at $s=1$. The parity of the analytic rank is
determined by the root number $w(E/K)\in\{\pm 1\}$, which is the conjectural
sign in the functional equation for $L(E/K,s)$ under $s\leftrightarrow 2-s$.
Although this $L$-function is not even known to exist at $s=1$
for $K\ne\Q$, the definition of the root number (due to Langlands)
is independent of any conjectures. Thus we expect the following parity
conjecture:

\begin{conjecture}
\label{parityconj}
For any (some) prime $p$, the root number agrees with the parity of
the $p$-Selmer rank, so $w(E/K)=\sigma(E/K,p)$.
\end{conjecture}

One of the main results in the paper is the following theorem.

\begin{theorem}
\label{main}
If $E/K$ has a rational isogeny of prime degree $p\ge 3$,
and $E$ is semistable at all primes over $p$,
then Conjecture \ref{parityconj} holds for $E/K$ and~$p$.
It also holds for $p=2$ under the additional assumption
that $E$ is not supersingular at primes above 2.
\end{theorem}

Recall that the global root number can be expressed in terms of
local root numbers over all places of $K$,
$$
  w(E/K)=\prod_v w(E/K_v) \>.
$$
If $E$ has an isogeny $\phi$ of degree $p$ over $K$,
then there is also a product formula for the parity of the $p$-Selmer rank
(Cassels-Fisher, see \cite{VD} Appendix),
$$
  \sigma(E/K,p) = \prod_v \sigma_\phi(E/K_v) \>.
$$
Here $\sigma_\phi(E/K_v)\in\{\pm 1\}$ is 1 if the power of $p$ in
$$
  \frac{\#\coker(\phi: E(K_v) \to E'(K_v))}{\#\ker(\phi: E(K_v) \to E'(K_v))}
$$
is even and $-1$ otherwise.

\begin{notation}
If $F$ is local, we write $(a,b)=(a,b)_F\in\{\pm 1\}$ for the Hilbert symbol:
it is 1 if and only if $b$ is a norm from $F(\sqrt{a})$ to $F$.

If $\phi:E\to E'$ is an isogeny defined over $F$,
we write $\f{F}$ for the extension of $F$ generated by the points
in $\ker\phi$. Since $\Gal(\f{F}/F)\hookrightarrow (\Z/p\Z)^*$
from the action on these points, $\f{F}/F$ is cyclic. We denote the image
of $-1$ under the composition
$$
  F^* \>\>{\buildrel\text{\tiny loc.recip.}\over
     {\hbox to 40pt{\rightarrowfill}}
  }\>\>\Gal(\f{F}/F)\hookrightarrow (\Z/p\Z)^*
$$
by $(-1,\f{F}/F)$, and we refer to it as the local Artin symbol.
It is 1 if $-1$ is a norm from $\f{F}$ to $F$ and $-1$ otherwise.
\end{notation}

In this paper we derive formulae for the local terms $w(E/K_v)$
and $\sigma_\phi(E/K_v)$ for odd $p$ (Theorems \ref{thm1} and \ref{thm2}).
It turns out that although the root number and the $p$-Selmer rank agree
globally, the local terms are not the same but are related as follows.

\begin{theorem}
\label{main2}
Let $K$ be a number field and $p$ an odd prime. Let $E/K$
be an elliptic curve with a cyclic $p$-isogeny $\phi$
defined over $K$, and assume that $E$ has semistable reduction
at all primes above $p$. Then for all places $v$ of $K$,
$$
  w(E/K_v) = (-1,\fKv/K_v)\>\sigma_\phi(E/K_v)\>.
$$
\end{theorem}

Note that Theorem \ref{main2} implies Theorem \ref{main}
by the product formula for the local Artin symbols, $\prod_v (-1,\fKv/K_v)=1$.

When $p=2$, the above theorem does not hold. The existence of a 2-isogeny
is equivalent to having a 2-torsion point, so the extension $\fKv/K_v$ is
always trivial. However, there is an analogous relation between
$w(E/K_v)$ and $\sigma_\phi(E/K_v)$ as follows.

Translating the 2-torsion point to $(0,0)$, the curves $E$ and $E'$
get the models
\begin{eqnarray}
  \label{abcurve}
  E:  & y^2 = x^3 + a x^2 + b x, \qquad a,b\in \O_K\>,  \\
  \label{abisocurve}
  E': & y^2 = x^3 - 2ax^2 + \delta x, \qquad \delta=a^2-4b\>,
\end{eqnarray}
with the isogeny $\phi: E \to E'$ given by
\begin{equation}
\label{abisogeny}
  \phi: (x,y) \mapsto (x+a+bx^{-1},y-bx^{-2}y) \>.
\end{equation}

\begin{theorem}
\label{p2thm}
Suppose $E/K$ has either good ordinary or multiplicative reduction at
all primes above 2. Then for all places $v$ of $K$,
\beql{2formula}
  w(E/K_v) = \sigma_\phi(E/K_v) \> (a,-b)_{K_v} \> (-2a,\delta)_{K_v} \>.
\eeql
In particular, the 2-parity conjecture holds for $E/K$ by the product formula
for the Hilbert symbols.
\end{theorem}


For $K=\Q$, the parity conjecture for $E$ and $p$ in the case that
$E$ has a $p$-isogeny is a theorem of P. Monsky \cite{Mon}, who also proved
it unconditionally for $K=\Q, p=2$.
J.~Nekov\'a\v r \cite{Nek} proved the conjecture
without the assumption on the existence of a $p$-isogeny
for elliptic curves over $\Q$
with potentially ordinary or potentially multiplicative reduction at $p$.
If $E/\Q$ is semistable and has a rational $p$-isogeny ($p$ odd),
the parity conjecture for $E$ base changed to an arbitrary number field
follows from \cite{VD}, Thm. 3 and Prop A.1.
Indeed, our computations of Selmer ranks are based on the approach by
T. Fisher in \cite{VD}. We would also like to mention that
for $E/\Q$, recently
M. Shuter \cite{Shu} has done some beautiful computations
of Selmer ranks over the fields where they acquire a $p$-isogeny.

\subsection*{Acknowledgements}
We would like to thank the referee for carefully reading the manuscript
and numerous comments and suggestions.

\section{Main results for $p \ge 3$}

\begin{theorem}
\label{thm1}
Assume $F=\R$ or $\C$, or $[F:\Q_l]<\infty$, and let $p\ge 3$.
Let $E/F$ be an elliptic curve with a rational $p$-isogeny
$\phi$.
Then
$$
  w(E/F)=\left\{\begin{array}{ll}
    -1 & \text{$F$ is Archimedean,}\cr
    1, & \text{$E$ has good reduction,}\cr
    -1, & \text{$E$ has split multiplicative reduction,}\cr
    1, & \text{$E$ has non-split multiplicative reduction,}\cr
    \delta\,\cdot\,(-1,\f{F}/F), & \text{$E$ has additive reduction and $l\ne p$.}\cr
  \end{array}\right.
$$
Here $\delta=1$ unless $p=3$, $\mu_3\not\subset F$ and $E/F$ has
reduction type IV or IV$^*$, in which case $\delta=-1$.
\end{theorem}

\begin{proof}
Except in the case of additive reduction, the formula for $w(E/F)$
is well-known and does not depend on the existence of a rational isogeny
(see e.g. \cite{RohG} Thm. 2).
The remaining case is dealt with in Section \ref{sRoot}.
\end{proof}

\begin{theorem}
\label{thm2}
Assume $F=\R$ or $\C$, or $[F:\Q_l]<\infty$, and let $p\ge 3$.
Let $E/F$ be an elliptic curve with a rational $p$-isogeny
$\phi$.
Then
$$
  \sigma_\phi(E/F)=\left\{\begin{array}{ll}
    -(-1,\f{F}/F), & \text{$F$ is Archimedean,}\cr
    \>\>\>(-1,\f{F}/F), & \text{$E$ has good reduction,}\cr
    -(-1,\f{F}/F), & \text{$E$ has split multiplicative reduction,}\cr
    \>\>\>(-1,\f{F}/F), & \text{$E$ has non-split multiplicative reduction,}\cr
    \delta, & \text{$E$ has additive reduction and $l\ne p$.}\cr
  \end{array}\right.
$$
Here $\delta=1$ unless $p=3$, $\mu_3\not\subset F$ and $E/F$ has
reduction type IV or IV$^*$, in which case $\delta=-1$.
\end{theorem}

\paragraph{{\bf Remark on Artin symbols}}
For $l\ne p$ the above local Artin symbols can be easily described:
If $F$ is Archimedean, $(-1,\f{F}/F)=1$ (i.e. $-1$ is a norm from $\f{F}$)
unless $F=\R$ and $\f{F}=\C$.
If $F$ is non-Archimedean and $E$ has semistable reduction,
then $(-1,\f{F}/F)=1$
because this extension is unramified (see proof below).
For $l=p$, see Lemma \ref{lemhilp} for the description of the Artin symbol.

\begin{proof}[Proof of Theorem \ref{thm2}]

Recall that $\sigma_\phi(E/F)=\pm 1$ and it is 1 if and only if the power of
$p$ in
\beq{kercoker}
  \frac{\#\coker(\phi: E(F) \to E'(F))}{\#\ker(\phi: E(F) \to E'(F))}
\eeq
is even.
For Archimedean $F$, the cokernel is always trivial while $\#\ker=p$
unless $F=\R$ and $\f{F}=\C$;
so $\sigma_\phi(E/F)=-(-1,\f{F}/F)$.

Henceforth assume that $F$ is a finite extension of $\Q_l$.
Then \eqref{kercoker} equals (\cite{VD}, App., formula (18))
$$
  \frac{c_v(E')}{c_v(E)}\cdot |\alpha_v|_v^{-1}
$$
where $c_v$ is the Tamagawa number and
$\alpha_v$ is the leading coefficient for the action of $\phi$ on the
formal groups. We will compute both contributions.

For the quotient ${c_v(E')}/{c_v(E)}$, Lemma \ref{lemtam}
in Section \ref{sTam} shows that it has odd $p$-valuation precisely
for the primes of split multiplicative reduction and primes
of additive reduction with $\delta=-1$.
Next, for $l\notdiv p$, the isogeny $\phi$
induces an isomorphism on formal groups, so $\alpha_v$ is a unit. For $l|p$,
we will show that $\ord_p|\alpha_v|_v$ is even if and only if
$(-1,\f{F}/F)=1$ (Section \ref{sAlpha}).

To complete the proof of the theorem, it remains to show that $(-1,\f{F}/F)=1$
for places $l\notdiv p$ of semistable reduction. It suffices to check
that $\f{F}/F$ is unramified, since then all units are norms by local
class field theory.
But $\f{F}/F$ is a Galois extension of degree prime to $p$,
while the inertia subgroup of $\Gal(K_v(E[p])/K_v)$
is either trivial in case of good reduction or a $p$-group
in case of multiplicative reduction (cf. \cite{Sil2}, Exc. 5.13).
\end{proof}

\begin{corollary}
\label{wloc}
Let $K$ be a number field and $E/K$ an elliptic curve with semistable
reduction at all primes above $p$. Assume $E$ has a rational $p$-isogeny
$\phi$.
Then
$$
  w(E/K) = \sigma(E/K,p) = (-1)^{\#\{v|\infty\}} (-1)^s
    \prod_{\text{$v$ additive}} \delta_v\>(-1,\fKv/K_v),
$$
where $s$ is the number of primes of split multiplicative reduction of $E/K$,
and $\delta_v=1$ unless $p=3$, $\mu_3\not\subset K_v$ and $E$ has
reduction type IV or IV$^*$ at $v$, in which case $\delta_v=-1$.
\end{corollary}

Since the local Artin symbol over the semistable primes $v\notdiv p$
is trivial, it follows from the product formula
that the product over the additive primes can be
replaced by the ones over $p$ and over $\infty$,
apart from the easy correction terms $\delta_v$ (trivial for $p>3$),
$$
  w(E/K) = \sigma(E/K,p) = (-1)^{\#\{v|\infty\}} (-1)^s
    \prod_{v|p\>\text{or}\>\infty} (-1,\fKv/K_v) \prod_{\text{$v$ additive}}\delta_v\>.
$$

\section{Root numbers}
\label{sRoot}

In this section $[F:\Q_l]<\infty$ and $E/F$ is an elliptic curve
with additive reduction which admits a cyclic $p$-isogeny
for some odd $p\ne l$. To prove Theorem \ref{thm1} we need to show that
$$
  w(E/F) = \delta\,\cdot\,(-1,\f{F}/F),
$$
where $\delta=1$ unless $p=3$, $\mu_3\not\subset F$ and $E/F$ has
reduction type IV or IV$^*$, in which case $\delta=-1$.
Recall also that $\f{F}$ is the extension of $F$ generated by the points
in the kernel of $\phi$.

We will determine $w(E/F)$ from the action of $\Gal(\bar F/F)$ on the
$p$-adic Tate module $T_p(E)$. We mention that computations of this kind
have previously been carried out by Rohrlich \cite{RohE, RohG} and
Kobayashi \cite{Kob}, and we refer to them for definitions
and background for local root numbers and $\epsilon$-factors
of elliptic curves.

We set $V_p(E)=T_p(E)\tensor_{\Z_p}{\bar\Q_p}$,
and recall
that the Weil group of $F$ is the subgroup of $\Gal(\bar F/F)$ generated
by the inertia subgroup and a lifting $\Frob$ of the Frobenius element.
Write $||\cdot||$ for the cyclotomic character (local reciprocity map
composed with the normalized absolute value of $F$).


\begin{lemma}
\label{lempotmult}
Suppose $E$ has potentially multiplicative reduction. Then
the action of the Weil group on $V_p(E)$ is of the form
$\smallmatrix{\chi}{*}{0}{\>\chi^{-1}||\cdot||}$
for some ramified quasi-character $\chi$, and inertia acts via
$\pm\smallmatrix1*01$
with $*$ not identically 0. The root number is given by $w(E/F)=(-1,\f{F}/F)$.
\end{lemma}

\begin{proof}
That inertia acts as asserted follows the theory of the Tate curve
(cf. \cite{Sil2} Lemma V.5.2, Excs. 5.11, 5.13).
In particular, $E$ acquires multiplicative reduction over $\f{F}$.

Because the inertia subgroup is normal in the Weil group, Frobenius
preserves the 1-dimensional subspace where inertia acts through a
quotient of order 2; this gives the action of the full Weil group.


Next, the root number of the semi-simplification of $V_p(E)$ is
given by the determinant formula (see \cite{RohE} p.145 or \cite{TatN} (3.4.7))
$$
  w(V_p(E)_{ss}) = w(\chi\oplus\chi^{-1}||\cdot||) = w(\chi\oplus\chi^{-1})
    = \chi(\theta(-1))
$$
with $\theta$ the local reciprocity map on $F^*$.
Over $\f{F}$ the quasi-character $\chi$, and therefore also $V_p(E)_{ss}$,
is unramified. Take a primitive
character $\tilde\chi$ of $\f{F}/F$ that coincides with $\chi$ on inertia.
Then
$$
  \chi(\theta(-1))=\tilde\chi(\theta(-1))=(-1,\f{F}/F) \>.
$$
The assertion follows from the formula (see \cite{TatN} (4.2.4))
\begin{equation}
  \epsilon(V_p(E)) = \epsilon(V_p(E)_{ss})
  \frac{\det(-\Frob^{-1}|{V_p(E)_{ss}}^{I})}{\det(-\Frob^{-1}|V_p(E)^{I})},
\end{equation}
since ${V_p(E)_{ss}}^{I}=0=V_p(E)^{I}$.
\end{proof}

\begin{lemma}
\label{lemadd5}
Suppose $E$ has potentially good reduction and $p\ge 5$. Then
the action of the Weil group on $V_p(E)$ is of the form
$\smallmatrix{\chi}{0}{0}{\>\chi^{-1}||\cdot||}$
for some quasi-character $\chi$.
The root number is given by $w(E/F)=(-1,\f{F}/F)$.
\end{lemma}

\begin{proof} From
the properties of the Weil pairing,
inertia acts via $\smallmatrix1*01$
on $E[p]$ over $\f{F}$.
On the other hand, the inertia is finite of order
dividing 24 (\cite{Ser} \S5.6), so has no elements of order $p$.
So it acts trivially on $E[p]$, hence $E/\f{F}$ has good reduction
(by \cite{ST} Cor. 2 or \cite{Sil2} Prop 10.3).

We need to show that the action of the Weil group on $V_p(E)$ is abelian.
On the one hand, the commutator of any two elements acts trivially on
the residue field, so it is an element of the inertia subgroup.
On the other hand,
its image in $\Gal(\f{F}/F)$ is trivial, because the latter is abelian.
As $E/\f{F}$ has good reduction, this commutator acts trivially on $T_p(E)$.

It follows that $T_p(E)\iso\chi\oplus\chi^{-1}||\cdot||$, so
$w(E/F)=(-1,\f{F}/F)$ as in the proof of Lemma \ref{lempotmult}.
\end{proof}

\begin{lemma}
Suppose $E$ has potentially good reduction and $p=3$. Then
$w(E/F)=\delta(-1,\f{F}/F)$.
\end{lemma}

\begin{proof}
Denote $G=\Gal(F(E[3])/F)$ and write $I$ for its inertia subgroup.
Since $I$ is a non-trivial subgroup of
$\smallmatrix**0*\subset\GL_2(\F_3)$ of
determinant 1, it is one of $C_2$, $C_3$ and $C_6$.
Moreover, $I=C_3$ if and only if $E$ has reduction type IV or IV$^*$.
(For $l\ne 2,$ $E$ has type IV or IV$^*$ if and only if the valuation
of the minimal discriminant of $E$ is 4 or 8, equivalently $|I|=3$.
For $l=2$, see \cite{Kra}, Thm. 2(i). )

(a) If $I=C_2$, then $\delta=1$. The root number is
$(-1,\f{F}/F)$ by the same argument as in Lemma \ref{lemadd5}.

(b) If $I=C_3$, then
$(-1,\f{F}/F)=1$ because it corresponds to an element of $I$
of order dividing 2. Next, $G$ is either $C_3$ or $C_6$
if $\mu_3\subset F$, and $S_3$ otherwise.

If $G$ is cyclic, then $w(E/F)=1$, because
$E$ acquires good reduction after a Galois cubic extension, and
the root number of an elliptic curve is unchanged in such an extension
(\cite{KT}, proof of Prop. 3.4).

If $G=S_3$, then $\Frob^2$ acts centrally on $V_3(E)$, so
it is given by
a scalar matrix $\lambda\,\text{Id}$. From
the properties of the Weil pairing, its determinant $\lambda^2$ is equal
to $\det(\Frob)^2=f^2$, where $f$ is the size of the residue field of $F$.
Note that $f\equiv 2\mod3$ as $\mu_3\not\subset F$ and
$\lambda\equiv1\mod 3$ since $\Frob^2$ acts
trivially on $E[3]$. In other words $\lambda=-f$.

Let $\chi$ be the unramified quasi-character of the Weil group
that takes Frobenius to $1/\sqrt{-f}\in\Q_3$. Then $V_3(E)\tensor\chi$
coincides with the 2-dimensional irreducible representation of $G\iso S_3$.
If $\phi$ is a character of order 3 of $I$, then
by a theorem of Fr\"ohlich-Queyrut (\cite{FQ}, Lemma 1 and Thm. 3),
$$
  w(\phi)=\phi(\theta(\sqrt{-3}))=1 \>,
$$
where $\theta$ is the local reciprocity map on $F(\mu_3)^*$.

Let $\eta$ be the quadratic unramified character of $G$, and
denote by $m$ the largest integer such that
$\tr_{F/\Q_l}(\pi_F^{-m}\O_F)\subset\Z_l$.
Writing $\1$ for the trivial representation,
by inductivity in degree 0,
$$
  1 = \frac{w(\phi)}{w(\1_{F(\mu_3)})}
    = \frac{w(V_3(E)\tensor\chi)}{w(\1_F)w(\eta)}
    = \frac{w(V_3(E)\tensor\chi)}{\eta(\Frob)^{m}} = (-1)^m w(V_3(E)\tensor\chi)\>.
$$
On the other hand, by the tensor product formula (\cite{TatN} (3.4.6)),
$$
  \epsilon(V_3(E)\tensor\chi)=\chi(\Frob^{-1})^{n(E)+2m}\epsilon(V_3(E))\>,
$$
so $w(E/F)=(-1)^{n(E)/2}=-1$ because $n(E)=2$ (tame ramification).


(c) Now assume that $I=C_6$. Then $G=C_6$ if $\mu_3\subset F$ and
$G\iso D_{12}$ otherwise. In the first case, the action of the Weil group
is abelian, so the same argument as in Lemma \ref{lemadd5} applies,
and $w(E/F)=(-1,\f{F}/F)$.

Finally, suppose $G\iso D_{12}$ and consider the twist $E_\chi$ of $E$
by the non-trivial character $\chi$ of the
quadratic extension $\f{F}/F$.
By inductivity in degree 0,
$$
  \frac{w(E/\f{F})}{w({\1_{\f{F}}\oplus\1_{\f{F}}})}
    = \frac{w(E/F)w(E_\chi/F)}{w(\1_F)^2 w(\chi)^2}\>.
$$
Both $E/\f{F}$ and $E_\chi/F$ have root number $-1$, as
they fall under case (b) with $G=S_3$. By the determinant formula,
$$
  w(E/F)=w(\chi)^2=w(\chi\oplus\chi^{-1})=\det(\chi)(-1)=(-1,\f{F}/F)\>.
$$
\end{proof}

\section{Tamagawa numbers}
\label{sTam}

In the next 3 sections, we complete the proof of Theorem \ref{thm2}.

\begin{lemma}
\label{lemtam}
Let $E/F$ be an elliptic curve, $[F:\Q_l]<\infty$. Suppose
$\phi: E\to E'$ is a cyclic $p$-isogeny defined over $F$ with $p\ge 3$.
Denote by $c(E)=[E(F):E_0(F)], c(E')=[E'(F):E'_0(F)]$ the Tamagawa numbers,
and let $\delta$ be as in Theorems \ref{thm1} and \ref{thm2}. Then
$$
  \ord_p \frac{c(E')}{c(E)} = \left\{
  \begin{array}{ll}
    0,      & \text{$E$ has good or non-split multiplicative reduction}\cr
    \pm1,   & \text{$E$ has split multiplicative reduction}\cr
    0,      & \text{$E$ has additive reduction and $\delta=1$}\cr
    \pm1,   & \text{$E$ has additive reduction and $\delta=-1$}
  \end{array}
  \right.
$$
\end{lemma}

\begin{proof}
If $E$ (and therefore also $E'$) has good reduction, then $c(E)=c(E')=1$.
If $E$ has non-split multiplicative reduction, then the $c$ are 1 or 2, so
the quotient is prime to $p$. If the reduction is split multiplicative,
the quotient contributes either $p$ or $p^{-1}$ (\cite{VD}, Lemma A.2).
If $E$ has additive reduction then $1\le c\le 4$,
so the quotient is prime to $p$ for $p\ge 5$.
It suffices to prove that for $p=3$
the quotient is prime to 3 precisely when $\delta=1$.

{}From Tate's algorithm (\cite{Sil2}, IV.9), the case $c=3$
only occurs when the reduction type is IV or IV$^*$.
Applying the multiplication-by-3 map to the exact sequence
$$
  0\lar E_0(F)\lar E(F)\lar E(F)/E_0(F)\lar 0
$$
and recalling that it is an
isomorphism on formal groups,
we get that $E(F)$ has a 3-torsion point if and only if
$c=3$.

If the absolute Galois group acts on $E[3]$ via
$\smallmatrix{\chi_1}{*}{0}{\chi_2}$, then its action on $E'[3]$ is of
the form $\smallmatrix{\chi_1}{0}{*}{\chi_2}$. Also note that
$\mu_3\subset F$ if and only if the action factors through $\SL_2(\F_3)$.
So if $\mu_3\subset F$, then $E(F)$ has a 3-torsion point if and only
if the isogenous curve has one, so $c(E)/c(E')=1$.
Conversely, if $\mu_3\not\subset F$, exactly one of $E(F), E'(F)$ has
a 3-torsion point, so $c(E)/c(E')=3^{\pm 1}$.
\end{proof}

\section{Local Artin symbols at primes above $p$}

\begin{lemma}
\label{lemhilp}
Let $\Q_p\subset F\subset F'$ be finite extensions ($p$ odd),
with $F'/F$ cyclic Galois of degree dividing $p-1$.
Then $(-1,F'/F)=1$ if and only if one of
the following conditions is satisfied:
\begin{enumerate}
\item The residue field $k_F$ of $F$ is of even degree over $\F_p$, or
\item $(p-1)/e(F'/F)$ is even, where $e$ denotes ramification degree.
\end{enumerate}
\end{lemma}

\begin{proof}
The condition $(-1,F'/F)=1$ is equivalent to $-1$ being a norm from $F'$ to $F$.
If $F_0$ is the maximal odd degree extension of $F$ inside $F'$, then
$N_{F_0/F}(-1)=-1$, implies $(-1,F'/F)=(-1,F'/F_0)$. In other words,
we may assume that $[F':F]$ is a power of 2.

Let $F^u$ be the maximal unramified extension of $F$ inside $F'$.
If $F^u=F'$, then all units of $F$ are norms from $F'$ and the result holds.
Otherwise, we can write $-1=\zeta^{[F^u:F]}=N_{F^u/F}(\zeta)$ for some
$\zeta\in\mu_{p-1}\subset F$. Then $(-1,F'/F)=1$ if and only if $\zeta$ is
a norm from $F'$ to $F^u$.

Since $F'/F^u$ is a totally and tamely ramified extension,
a unit in $F^u$ is a norm from $F'$ if and only if its reduction lies in the
unique subgroup of $k^*$ of index $[F':F^u]$, where $k$ is the residue field
of $F^u$. 

Writing $d=[F^u:F]$, we have
\begin{equation}
\label{ineq1}
  \ord_2\Bigl(\frac{p-1}{2d}\Bigr) = \ord_2 [\F_p^* : \langle\bar\zeta\rangle]
    \le \ord_2 [k^* : \langle\bar\zeta\rangle]\>,
\end{equation}
where the last inequality is an equality if and only if $[k:\F_p]$ is odd,
equivalently $d=1$ and $[k_F:\F_p]$ is odd. Also,
\begin{equation}
\label{ineq2}
  \ord_2\Bigl(\frac{[F':F]}{d}\Bigr) \le \ord_2\Bigl(\frac{p-1}{d}\Bigr)
    = \ord_2\Bigl(\frac{p-1}{2d}\Bigr) + 1\>,
\end{equation}
the first equality holding if and only if $(p-1)/[F':F]$ is odd.
On the other hand,
$$
  \langle\zeta\rangle \not\subset N_{F'/F^u} {F'}^*  \quad\iff\quad
  \ord_2[k^*:\langle\bar\zeta\rangle] < \ord_2[F':F^u] =
    \ord_2\Bigl(\frac{[F':F]}{d}\Bigr)  \>.
$$

If both the conditions (1) and (2) in the lemma fail,
then $(p-1)/[F':F]$ is odd, $F'/F$ is totally ramified (so $d=1$),
and the inequalities in \eqref{ineq1}, \eqref{ineq2} become equalities.
Hence $\langle\zeta\rangle \not\subset N_{F'/F^u} {F'}^*$ and $(-1,F'/F)\ne 1$.

Conversely, if one of (1) and (2) is satisfied, one of
the inequalities in \eqref{ineq1}, \eqref{ineq2} is strict, so
$\ord_2[k^*:\langle\bar\zeta\rangle] \ge \ord_2[F':F^u]$. In other words,
$\zeta$ is norm from $F'$ and $(-1,F'/F)=1$.
\end{proof}

\section{$p$-isogenies on formal groups}
\label{sAlpha}

Let $F$ be a finite extension of $\Q_p$ ($p$ odd),
and denote by $\O_F, \m=(\pi_F), v$ and $k_F=\O_F/\m$
its ring of integers, the maximal ideal, the valuation and the residue field respectively.
Suppose $E/F$ is an elliptic curve with semistable
reduction and let $\phi: E\to E'$ be a cyclic $p$-isogeny defined over $F$.
Let $f: \cF_E(\m) \to \cF_{E'}(\m)$ be the induced map on the formal groups,
which can be considered as a power series of the form
$$
  f(T) = \alpha T + \ldots \>.
$$
Write $|\alpha|_v=p^{-[k_F:\F_p]v(\alpha)}$
for the normalized absolute value of $\alpha$ in $F$ and let
$2^e||p-1$.
We claim that $|\alpha|_v$ is an odd power of $p$
if and only if $2^e$ divides the ramification index of $\f{F}/F$ and
$[k_F:\F_p]$ is odd. By Lemma \ref{lemhilp}, this will complete
the proof of Theorem~\ref{thm2}.


First of all, if $[k_F:\F_p]$ is even
then clearly $|\alpha|_v$ is an even power of $p$ and the statement
holds. So suppose now that $[k_F:\F_p]$ is odd, in which case
$\ord_p |\alpha|_v \equiv v(\alpha)\mod 2$.


Let $\bar f: \overline{\cF_E} \to \overline{\cF_{E'}}$
be the reduction of $f$ modulo $\m$.
The reduced formal groups $\overline{\cF_E}, \overline{\cF_{E'}}$
either are those of the reduced elliptic curve in the case
of good reduction, or become isomorphic to ${\mathbb G}_m$ over $F^{nr}$
in the case of multiplicative reduction
(as follows from the theory of the Tate curve, cf. \cite{Ser}, p.277).
The map $\bar f$ is an isogeny of formal groups over $k_F$
of degree dividing $p$, and it either has height 0 or 1
(see \cite{Sil1}, Thm. IV.7.4). We have two cases to consider:

\begin{lemma}
If $\alpha$ is a unit, then $\f{F}/F$ is unramified.
\end{lemma}

\begin{proof}
That $\alpha$ is a unit means that $\bar f$ is an isomorphism of the reduced
formal groups. Then the group scheme $\ker\phi$ is \'etale over $\O_F$,
so $\f{F}/F$ is unramified.
\end{proof}

In the case that $v(\alpha)>0$, the reduction $\bar f$ is an inseparable
isogeny of degree $p$. Then we have

\begin{lemma}
\label{lemfor}
If $\bar f$ is inseparable of degree $p$, then $f$ has a kernel
of order $p$ in the maximal unramified extension $F^{nr}$
if and only if $v(\alpha)$ is a multiple of $p-1$.
\end{lemma}

\begin{proof}
Let $\omega(T), \omega'(T)$ be the normalized invariant differentials on
$\cF_E, \cF_{E'}$. Then (\cite{Sil1}, Cor. IV.4.3)
$$
  \omega'\compose f = \alpha\,\omega \>,
$$
so $\alpha\omega(T)=(1+\cdots)\frac{d}{dT} f(T)\>.$
Because $(1+\cdots)$ is invertible, it follows that
$$
  f(T)=\alpha f_1(T) + f_2(T^p)
$$
for some $f_1, f_2 \in \O_F[[T]]$. (This is the same argument
as in \cite{Sil1}, Cor. IV.4.4.) Moreover $\bar f$ has height 1, so
$$
  f_1(T)=T+\ldots, \qquad f_2(T)=u T+\ldots \>\>\> (u\in\O_F^*)\>.
$$
Now we can prove the lemma.

$\Rightarrow$. We have a kernel of order $p$, so $f(m)=0$ for some
$m\in\mfnr, m\ne 0$. The first terms in the expansions for $f_1(m)$ and
$f_2(m^p)$ must cancel modulo $\mfnr$.
Hence $v(\alpha m)=v(um^p)$, so $v(\alpha)$ is divisible by $p-1$.

$\Leftarrow$. If $\alpha$ has valuation divisible by $p-1$, write
$\alpha=\alpha_0^{p-1}$ with $\alpha_0\in\mfnr$. Replace $T$ by $T\alpha_0$, so
$$
  f(\alpha_0T) = (\alpha_0)^{p-1} f_1(\alpha_0T) + f_2((\alpha_0T)^p)
       = \alpha_0^p ( T + ... + u T^p + ... ) = \alpha_0^p g(T),
$$
with every coefficient in $...$ having positive valuation.
Because $g'(T)$ is a unit, by Hensel's lemma the $p$ distinct roots
of $g(T)\mod\mfnr$ lift to roots $z_i$ of $g(T)$. Then $z_i\alpha_0$
are $p$ distinct roots of $f(T)$, so $f$ has a kernel of order~$p$.
\end{proof}

Now we can complete the proof of Theorem \ref{thm2}.
Because $\phi$ has $p$ points in the kernel over $\f{F}$, by the above lemmas
$2^e|v_{\f{F}}(\alpha)$. If $v(\alpha)$ is odd, this means that the ramification
degree of $\f{F}/F$ is a multiple of $2^e$, as asserted.
If $v(\alpha)$ is even, then over
$F^{nr}(\sqrt[\frac{p-1}{2}]{\pi_F})$
the map $f$ acquires a kernel (Lemma \ref{lemfor} again), so
$\f{F}\subset F^{nr}(\sqrt[\frac{p-1}{2}]{\pi_F})$
has ramification degree over $F$ not divisible by $2^e$.

\def\beq{$$\begin{array}{llllll}}
\def\eeq{\end{array}$$}

\section{Parity conjecture for curves with a 2-isogeny}

The purpose of this section is to prove Theorem \ref{p2thm}.
First note that $a\ne 0$, for otherwise $j(E)=1728$ and
$E$ has either additive or supersingular reduction at places above 2.

Henceforth we work in the local setting. Let
$F=\R$ or $\C$, or $[F:\Q_l]<\infty$.
In the non-Archimedean case we write $\m$ for the maximal ideal and
$v$ for the normalized valuation on $F$.
Suppose $E,E'/F$ are elliptic curves given by equations \eqref{abcurve} and
\eqref{abisocurve} with $a,b\in\O_F$ and $a\ne 0$.
Let $\phi: E \to E'$ be the 2-isogeny with kernel $\O,(0,0)$ and
defined by \eqref{abisogeny}.
The discriminant
$$
  \Delta(E) = 16\,\delta\,b^2
$$
is non-zero, so $b$ and $\delta$ are non-zero as well.
In particular, the Hilbert symbols $(a,b)$ and $(-2a,\delta)$
make sense, and we need to prove that
$$
  w(E/F) = \sigma_\phi(E/F) \> (a,-b) \> (-2a,\delta) \>.
$$

\subsection{Infinite places}

First suppose that $F=\R$ or $\C$, so $w(E/F)=-1$, and $\sigma_\phi(E/F)=1$
if and only if $\ord_2(\#\ker\phi/\#\coker\phi)$
is even. Clearly we need to show that
$$
  (a,-b) \> (-2a,\delta) =1 \quad\iff\quad
    \phi: E(F)\to E'(F) \text{ surjective.}
$$
If $F=\C$, then the Hilbert symbols are trivial and $\phi$ is surjective.

Suppose $F=\R$. If $(-2a)^2-4\delta=16b<0$,
then $E'(\R)\iso S^1$, so
$\phi$ is always surjective. On the other hand $-b>0$ implies
$(a,-b)=1$, and $\delta=a^2-4b>0$ implies $(-2a,\delta)=1$.

Similarly, if $b>0$ and $\delta<0$ then $E(\R)\iso S^1$ and
$E'(\R)\iso S^1\times\Z/2\Z$, so $\phi$ is not surjective; here
exactly one of the Hilbert symbols is 1, depending on the sign of $a$.

Finally, if $b, \delta>0$, then $E(\R)\iso S^1\times\Z/2\Z\iso E'(\R)$
and $(-2a,\delta)=1$.
Here $\phi$ is surjective if and only
if the points $\O,(0,0)$ of $\ker\phi$ lie on the same connected component.
(If they are on different components, the image of $\phi$ is connected;
otherwise, the identity component of $E(\R)$ maps 2-to-1 to the identity
component of $E'(\R)$, so the other component maps
to the other component since $\deg\phi=2$.)
So $\phi$ is surjective if and only if $0$ is the rightmost root of
$x^3+ax^2+bx$. This is equivalent to $-a<0$ and hence to $(a,-b)=1$.

\subsection{Finite places}

From now suppose that $F$ is a finite extension of $\Q_l$.
The only transformations that preserve the chosen model for $E$ are
$(x,y)\mapsto (u^4 x,u^6y)$. The constituents in the Hilbert
symbols then get multiplied by squares ($u^2,u^4$), and
the Hilbert symbols are unchanged.
So for $l\ne 2$ (including $l=3$)
the model \eqref{abcurve} may and will be chosen
to be minimal for the proof.

To prove the theorem, we will proceed as in the proof of Theorem \ref{thm2}.
Recall that $\sigma_\phi(E/F)=\pm 1$ and it is 1 if and only if the power of
$2$ in
$$
  \frac{\#\coker(\phi: E(F) \to E'(F))}{\#\ker(\phi: E(F) \to E'(F))}
$$
is even. The quotient equals 
$\tfrac{c(E')}{c(E)}\cdot|\alpha|_F^{-1}$ (\cite{VD}, App., formula (18))
where $c$ is the Tamagawa number,
$\alpha$ is the leading coefficient for the action of $\phi$ on the
formal groups, and $|\cdot|_F$ is the normalized absolute value.
We will compute $c(E)$, $c(E')$, $w(E/F)$ and $|\alpha|_F$.
For the latter, $\phi^*(dx/y)=dx/y$
from the explicit formula \eqref{abisogeny} for $\phi$.
So if $E, E'$ are transformed to their
respective minimal models by standard substitutions
$(x,y)\to (w^2x+...,w^3y+...)$ and $(x,y)\to (u^2x+...,u^3y+...)$, then
$\alpha=uw^{-1}$. We distinguish between various possibilities for the
reduction types:

\subsection{Good reduction, $l\ne 2$}

Here $w(E/F)=1$, $\sigma_\phi(E/F)=1$ and $b,\delta\in\O_{F}^*$.
If $a\in\O_{F}^*$, then both the Hilbert symbols are $($unit,unit$)$,
hence trivial. For $a\equiv0\mod\m$, the expression
$-b\delta\equiv 4b^2\mod\m$ is a non-zero square mod $\m$,
so the product of the Hilbert symbols is again trivial.

\subsection{Additive reduction, $l\ne 2$}

Reduction is either potentially multiplicative or potentially good.
In the latter case, $E/F(E[4])$ has good reduction and $F(E[4])/F$ is a
2-extension, so $3|v(\Delta)$.
\def\pl{^+}
We have the following options (writing $n\pl$ for an integer $\ge n$):

\smallskip
\begin{tabular}{lcccccc}
Reduction      & III & III${}^*$ & I${}_0^*$ & I${}_n^*\>\,(n>0)$ \cr
$w(E/F)$     & $(-2,\pi)$ & $(-2,\pi)$ & $(-1,\pi)$ & $(-1,\pi)$ \cr
$v(\Delta_E)$  & 3 & 9 & 6 & 6+n \cr
$v(a),v(b),v(\delta)\>$ & $\>1\pl,1,1\>$ & $\>2\pl,3,3\>$ & $\>1\pl,2,2\>$ & $\>1,2,3\pl$ or $1,3\pl,2$ \cr
\end{tabular}
\smallskip

\noindent
Here $\pi$ is a uniformiser of $F$, the reduction types are
from \cite{Sil2}, IV.9 and the root numbers are from \cite{RohG}
and \cite{Kob}. Note also that $E$ has potentially multiplicative reduction
(i.e. I${}_n^*$) if and only if $E'$ has; the same holds for I${}_0^*$
(inertia of order 2). In what follows we continue referring to the description
of Tate's algorithm in \cite{Sil2}, IV.9 for the description of the Tamagawa
numbers.

\paragraph{(III)} Here $c(E)=c(E')=2$. Next, $\delta\equiv -4b\mod\pi^2$, so
$$
  (a,-b)(-2a,\delta)=(a,-b\delta)(-2,\delta) = (a,4b^2+O(\pi^3))(-2,\pi)
    = (-2,\pi).
$$
\paragraph{(III$^*$)} Same argument, replacing $\pi^2$ and $O(\pi^3)$ by
$\pi^4$ and $O(\pi^5)$ respectively.

\paragraph{(I$_0^*$)} By assumption,
$T^3+\frac{a}{\pi}T^2+\frac{b}{\pi^2}T$ has
3 distinct roots mod $\m$ over the algebraic closure, and $b=\pi^2u, \delta=\pi^2w$
with $u,w$ units. Moreover $w$ is a square if and only if all three roots are
defined over the residue field of $F$, equivalently $c(E)=4$ (otherwise
$c(E)=2$). Similarly, $c(E')=4$ if and only if $u$ is a square, and 2
otherwise. So
\beql{addred1}
\begin{array}{lllll}
  \ord_2 c(E)\equiv 0 \mod 2&\iff&
  w \in (K^*)^2 &\iff& (\pi,\delta)=1, \cr
  \ord_2 c(E')\equiv 0 \mod 2&\iff&
  u \in (K^*)^2 &\iff& (\pi,b)=1.
\end{array}
\eeql
Now,
\beql{addred2}
  (a,-b)(-2a,\delta)=(a,-b)(-2,\pi^2w)(a,\delta)=(a,-b\delta)=(\pi,-b\delta)\>.
\eeql
The last equality is clear if $v(a)$ is odd; for $v(a)\ge 2$ even
$(a,-b\delta)=1$, and
$$
  -b\delta = -b(a^2-4b) \equiv 4b^2 \mod \pi^5
$$
is a square, so the last Hilbert symbol is 1 as well. Combining
\eqref{addred1}, \eqref{addred2} with $w(E/F)=(\pi,-1)$ yields the result.

\paragraph{(I$_n^*$)} We have $v(a)=1, v(b)\ge 2$ and so $v(\delta)\ge 2$.
Because $v(\Delta)=v(b^2\delta)>6$, one of $v(b), v(\delta)$ is at least 3,
hence the other one is 2, as $\delta+4b=a^2$; so we have essentially two cases.
Swapping $E$ and $E'$ interchanges $b$ and $\delta$ up to units,
so $E$ and $E'$ will always be in two different cases. We begin by
determining $c(E)\in\{2,4\}$.

Suppose $v(b)>2$ and $v(\delta)=2$, so
$E$ has type I${}_{2n}^*$ with $n=v(b)-2$.
Following \cite{Sil2}, IV.9, Step 7, the reduction of the polynomial
$P(T)=T^3+\frac{a}{\pi}T^2+\frac{b}{\pi^2}T$ mod $\m$ has a double root at
the origin; furthermore, $\frac{a}{\pi}X^2+\frac{b}{\pi^{2+n}}X$
has two distinct roots $\mod\m$, so $c(E)=4$.

Suppose $v(b)=2$ and $v(\delta)>2$, so $E$ has type I${}_{n}^*$
with $n=v(\delta)-2$. Translate $x$ by $a/2$ to get a model
$$
  y^2 = (x-\frac{a}{2})(x^2-\frac{\delta}{4}) \>,
$$
so that $P(T)=(T-\frac{a}{2\pi})(T^2-\frac{\delta}{4\pi^2})$
again has a double root at the origin mod $\m$.
Now, by the criterion in \cite{Sil2}, IV.9, Step 7, we have
$c(E)=4\Leftrightarrow\delta=\square$ for $n$ even, and
$c(E)=4\Leftrightarrow\frac{a}2\delta=\square$ for $n$ odd.

We are now in position to compute the Hilbert symbols and
to complete the proof in the I$_n^*$ case. To simplify the argument slightly,
note that showing \eqref{2formula} for $E$ is equivalent to that for $E'$:
the products of the Hilbert symbols
differ by $(-1,-2)=1$, both root numbers are $(\pi,-1)$ and
$\ord_2(c(E)/c(E'))\equiv\ord_2(c(E')/c(E))\mod 2$.
So we may assume without loss of generality that $v(b)<v(\delta)$, thus
$c(E')=4$. Note that in this case
$$
  b=a^2-4\delta=a^2(1-O(\pi))=\square.
$$
If $v(\delta)$ is even, then the parity of $\ord_2(c(E)/c(E'))$
is determined by the Hilbert symbol $(\pi,\delta)$, and
$$
  (a,-b)(-2a,\delta)=(a,-1)(-2,\delta)(a,\delta)=
    (\pi,-1)(\pi,\delta).
$$
Similarly, if $v(\delta)$ is odd, then the parity of $\ord_2(c(E)/c(E'))$
is determined by the Hilbert symbol $(\pi,2a\delta)$, and
$$
  (a,-b)(-2a,\delta)=(a,-1)(-2a,2a\delta)=(\pi,-1)(\pi,2a\delta).
$$

\subsection{A lemma on Hilbert symbols}

\begin{lemma}
\label{hil2}
Let $F/\Q_p$ be a finite extension. Then
\begin{enumerate}
\item $(1+4x,y)=1$ if $v(x)>0$ and $y\in F^*$,
\item $(1+4x,y)=1$ if $p=2$, $v(x)=0$ and $y\in\O_{F}^*$,
\item $(-1,-2)=-1$ if and only if $p=2$ and $[F:\Q_2]$ is odd.
\end{enumerate}
\end{lemma}

\begin{proof}
These statements are clear for odd $p$, so suppose that $p=2$.

1. \cite{Sil2}, Chapter V, Lemma 5.3.1.

2. It suffices to show that the extension $F(\sqrt{1+4x})/F$
is unramified, so every unit $y$ is a norm. Equivalently, if $L/F$ is the
unique quadratic unramified extension, we claim that $1+4x\in(L^*)^2$.
Let $\bar x\in\F_{2^n}$ be the reduction of $x$ mod $\m$.
Because
every quadratic polynomial over $\F_{2^n}$ has a root in $\F_{2^{2n}}$,
there is a unit $z$ of $L$ with $z^2+z\equiv x\mod \m$. Then
$$
  (1+2z)^2 = 1 + 4(z+z^2) \equiv 1+4x \mod 4\m,
$$
so $(1+4x)/(1+2z)^2$ is a square in $L$ by part 1, and so is $1+4x$.

3. If $\sqrt{-2}\in F$, then both of the conditions hold. Otherwise,
\beq
  (-1,-2)_K &=& (-1,F(\sqrt{-2})/F) = (N_{F/\Q_2}(-1),\Q_2(\sqrt{-2})/\Q_2)\cr
  &=& ((-1,-2)_{\Q_2})^{[F:\Q_2]} = (-1)^{[F:\Q_2]},
\eeq
as asserted.
\end{proof}

\subsection{Good reduction, $l=2$,
2-torsion point reduces to $\bar P=\O$}

Suppose $l=2$ and $E/F$ has good reduction,
so $w(E/F)=1$ and $c(E)=c(E')=1$. By assumption, the reduction
is ordinary, equivalently $j(E)$ is a unit.
By \cite{Sil1} A.1.1c, we can choose a minimal model of $E$ over $\O_{F}$
of the form
$$
  y^2+a_1xy+a_3y=x^3+a_2x^2+a_4x+a_6, \qquad 1\!-\!a_1,a_3,a_4\in\m\>.
$$
After the substitution $(x,y)\mapsto (a_1^2x-a_3a_1^{-1},a_1^3y)$,
we may assume $a_1=1$ and $a_3=0$.
Next, the substitution $(x,y)\mapsto (x-2t,y+t)$
with $t\in\m$ a root of $12t^2-(1+4a_2)t+a_4=0$
(Hensel's lemma) eliminates $a_4$.
Neither substitution changes the reduced curve,
thus we may assume that our model is
$$
  y^2+xy=x^3+a_2x^2+a_6, \qquad a_2\in\O_F,\>a_6\in\O_{F}^*.
$$
After completing the square, this becomes
$$
  y^2=x^3+(a_2+\frac14)x^2+a_6\>.
$$
Now let $(x_0,0)$ be our 2-torsion point with $v(x_0)<0$. Then
$$
  x_0^2(x_0+a_2+\frac14) = -a_6\in\O_{F}^*\>,
$$
so $v(x_0+a_2+1/4)=-2v(x_0)>0$, hence $v(x_0)=-v(4)$. Then
$v(x_0+a_2+1/4)=v(16)$, and write $x_0+a_2+1/4=16v$ with $v$ a unit.
Letting $w=1+4a_2$ and translating $x$ by $x_0$, the curve becomes
$$
  y^2 = x^3 + ax^2 + bx,
$$
with
\beq
  a=\frac12(-w+96v)=-\frac w2\cdot\square, \quad\cr
  b=\frac1{16}(-w+64v)(-w+192v)=\square, \quad\cr
  \delta=a^2-4b=-16v(-w+48w)=vw\cdot\square.
\eeq
Therefore
$$
  (a,-b)(-2a,\delta) = (-\frac w2,-1)(w,vw) = (-2,-1)(w,-vw) = (-2,-1) \>,
$$
where the last equality holds since (1 mod 4,unit)=1. On the other hand, the
isogenous curve
$$
  E': y^2 = x^3 - 2a x^2 + \delta x
$$
transforms under $x\to 4x, y\to 4x+8y$ to
$$
  y^2 + xy = x^3 + (a_2-24v)x^2 + (4a_2-48v+1) x,
$$
which has good reduction at 2. So $\ord_2|\alpha|_F=\ord_2|2|_F$ is even
if and only $[F:\Q_2]$ is even if and only if $(-2,-1)=1$
(Lemma \ref{hil2}).

\subsection{Good reduction, $l=2$, 2-torsion point reduces to $\bar P\ne\O$}

As before, $w(E/F)=1$ and $c(E)=c(E')=1$.
We claim that $\alpha=1$,
and that both Hilbert symbols are trivial.
Translating the 2-torsion point on the reduction to $(0,0)$, we
get a model
\beq
  E:  &y^2 + xy = x^3 + a_2 x^2 + a_4x, & a_2\in\O_{F},\>a_4\in\O_{F}^*,\cr
  E': &y^2 + xy = x^3 + a_2 x^2 - 4a_4x - (4a_2+1)a_4.
\eeq
These transform to our models \eqref{abcurve}, \eqref{abisocurve} with
substitutions $x\to x+..., y\to y+...$, so $\alpha=1$.
Next (cf. Lemma \ref{hil2} part 2),
\beq
  a&=&a_2+1/4,\cr
  b&=&a_4, \cr
  \delta&=&(a_2+1/4)^2-4a_4=1/16+a_2/2+a_2^2-4a_4, \cr
  (a,-b)&=&(a_2+1/4,-a_4)=(1+4a_2,-a_4)=(1\mod 4,\text{unit})&=1, \cr
  (-2a,\delta)&=&(-2a,1+8a_2+16a_2^2-64a_4)=(-2a,\square)&=1. \cr
\eeq

\subsection{Split multiplicative primes}

Write $E$ as a Tate curve (\cite{Sil2} \S V.3)
$$
  E_q: y^2+xy = x^3+a_4(q)x+a_6(q), \qquad E(F)\iso F^*/q^\Z,
$$
with $q\in m_K$ of valuation $v(q)=v(\Delta)$.  The coefficients have
expansions
$$
  a_4(q)=-5s_3(q), \quad a_6(q)=-\frac{5s_3(q)+7s_5(q)}{12}, \quad
    s_k(q)=\sum_{n\ge 1}\frac{n^k q^n}{1-q^n}\>,
$$
and they start
\beq
a_4(q)=-5q-45q^2-140q^3-365q^4-630q^5+O(q^6),\cr
a_6(q)=-q-23q^2-154q^3-647q^4-1876q^5+O(q^6)\>.
\eeq
The two-torsion, as a Galois set, is $\{1,-1,\sqrt{q},-\sqrt{q}\}$. For $u\ne 1$
in this set, the corresponding point on $E$ has coordinates
\beq
  X(u,q) = \frac{u}{(1-u)^2} + \sum\limits_{n\ge 1} \bigl(
    \frac{q^nu}{(1-q^nu)^2} + \frac{q^nu^{-1}}{(1-q^nu^{-1})^2} -2\frac{q^n}{(1-q^n)^2}
  \bigr),\cr
  Y(u,q) = \frac{u^2}{(1-u)^3} + \sum\limits_{n\ge 1} \bigl(
    \frac{q^{2n}u^2}{(1-q^nu)^3} + \frac{q^nu^{-1}}{(1-q^nu^{-1})^2} +\frac{q^n}{(1-q^n)^2}
  \bigr).
\eeq
We now have two cases to consider: the 2-torsion point
$(X(-1,q),Y(-1,q))\in E_q$ and (renaming $\pm\sqrt{q}$ by $q$)
the 2-torsion point $(X(q,q^2),Y(q,q^2))\in E_{q^2}$. In both cases,
we have $c(E)/c(E')=2^{\pm 1}$ and $w(E/F)=-1$, so we need
\beql{splitfor}
  \text{$\ord_2 |\alpha|_F$ even} \quad\iff\quad
    (a,-b) \> (-2a,\delta) = 1\>,
\eeql
where $a, b, \delta$ are the invariants of the curve \eqref{abcurve}, when
our curve is transformed into that from with the 2-torsion point at $(0,0)$.
First of all, $E_q$ has a model
$$
  y^2 = x^3+x^2/4+a_4(q)x+a_6(q)\>.
$$
Let $r=-X(u,q)$, and write $a_4=a_4(q), a_6=a_6(q)$. Then, after
translation, the curve becomes
$$
  E: y^2 = x^3+ax^2+bx, \quad a=1/4-3r, \quad b=2a_4-r/2+3r^2.
$$
Recall that the isogenous curve $E'$ is \eqref{abisocurve} in this form.

Suppose we are in Case 1, so $r=-X(-1,q)$. Then the substitution
\beql{mulsub}
  x\to 4x-2r+1/2, \quad y\to 8y+4x
\eeql
transforms $E'$ into the form
$$
  E^\dagger: y^2+xy = x^3 + (-5q^2+O(q^4)) x + (-q^2+O(q^4))\>.
$$
We use the notation $O(q^n)$ to indicate a power series in $q$ with
coefficients in $\O_{F}$ that begins with $a_nq^n+...$.
In fact, $E^\dagger=E_{q^2}$ but we will not need this; it is only
important that it is again a Tate curve (in particular, this model is minimal),
and $\alpha=2$ (this comes from \eqref{mulsub}). So
$$
  \text{$\ord_2|\alpha|_F$ even} \iff
  \text{$l\ne 2$ or $[F:\Q_2]$ is even} \iff (-1,-2)=1.
$$
Finally, from the expansions
$$
  r = 1/4+4 O(q), \quad
  a = -1/2+4 O(q), \quad
  b = 1/16+O(q),
$$
we have
\beq
  (a,-b) & = & (a,-1)(a,b)=(a,-1)(a,\square)=(-1/2+4O(q),-1) \cr
         & = & (-1/2,-1)(1+8O(q),-1) =(-2,-1)(\square,-1) = (-1,-2),\cr
  (-2a,\delta) & = & (1-8O(q),\delta) = (\square,\delta) = 1\>.
\eeq

Case 2 is similar and we omit the details;
here $\alpha$ is a unit, so we need to show that the
product of the two Hilbert symbols is 1. Here
$$
  a = 1/4+2\, O(q), \quad
  b = q+O(q^2), \quad
  \delta = 1/16 + O(q)\>.
$$
In particular, $a$ and $\delta$ are squares in $F$, so both Hilbert symbols
are trivial.

\subsection{Nonsplit multiplicative primes}

Let $F(\eta)/F$ be
the quadratic unramified extension of $F$ and consider
the twist of $E$ by $\eta$,
\beq
  E:      & y^2=x^3+ax^2+bx \>, \cr
  E_\eta: & y^2=x^3+\eta ax^2+\eta^2 bx \>.
\eeq
Then $E_\eta$ has split multiplicative reduction, so (cf. \eqref{splitfor})
$$
  \text{$\ord_2|\alpha_{E_\eta}|_F$ even} \iff (\eta a,-\eta^2 b) (2\eta a,\eta^2 (a^2-4b)) = 1.
$$
Also $\alpha_{E_\eta}=\alpha_E$, since the two curves become
isomorphic over $K_\eta$. Now,
\beq
  (\eta a,-\eta^2 b) = (\eta a,-b) = (\eta,-b) (a,-b), \cr
  (2\eta a,\eta^2 (a^2-4b)) = (2\eta a,a^2-4b) = (\eta,a^2-4b)(2a,a^2-4b),
\eeq
so comparing with the Hilbert symbols in \eqref{2formula} we have an extra term
\beql{2smcorr}
  (\eta,-b(a^2-4b)) = (\eta,-b^2\Delta(E')/16\Delta(E)) = (\eta,-\Delta(E')/\Delta(E)) \>.
\eeql
Because $x$ is a norm from $F(\eta)^*$ to $F^*$ if and only if $v(x)$ is even,
this Hilbert symbol is trivial precisely when
$v(\Delta(E'))\equiv v(\Delta(E))\mod 2$. From Tate's algorithm
(\cite{Sil2}, IV.9.4, Step 2),
$$
  c(E) = \left\{ \begin{array}{ll}
    1, & \text{$v(\Delta(E))$ odd},\cr
    2, & \text{$v(\Delta(E))$ even},
  \end{array}\right.
  \quad
  c(E') = \left\{ \begin{array}{ll}
    1, & \text{$v(\Delta(E'))$ odd},\cr
    2, & \text{$v(\Delta(E'))$ even},
  \end{array}\right.
$$
so the correction term \eqref{2smcorr} is trivial if and only
$c(E)/c(E')$ has even 2-valuation. This proves \eqref{2formula} in the
non-split multiplicative case.

\end{document}